\documentclass[11pt]{article}
\usepackage{amsmath}
\numberwithin{equation}{section}
\usepackage{amsfonts,amssymb,latexsym,amsbsy, bbm, theorem,enumerate,color}
\usepackage{amssymb}
\usepackage{graphicx, graphics}
\usepackage{float,color,fancybox,shapepar,setspace,hyperref}
\usepackage{subfigure}
\usepackage{pgf,tikz}
\usetikzlibrary{arrows}
\makeatletter

\newcommand{\Rmnum}[1]{\expandafter\@slowromancap\romannumeral #1@}
\makeatother
\newtheorem{Main Theorem}{Main Theorem}

\newtheorem{Theorem}{Theorem}

\newtheorem{Lemma}{Lemma}
\newtheorem{Claim}{Claim}

\def\square{\hbox{\vrule height8pt depth0pt
\vbox{\hrule width7.2pt\vskip7.2pt\hrule width7.2pt}\vrule
height8pt depth0pt}\smallskip}

\def\pf{\medskip\noindent {\emph{\bf Proof}.}~~}

\def\mytextindent#1{\indent\llap{#1\enspace}\ignorespaces}

\def\myitem{\par\hangindent\parindent\mytextindent}

\newcommand{\ex}{{\rm  ex}}

\usepackage{authblk}

\begin{document}

\title{Maximum cliques in a graph without disjoint given subgraph}
\author[1]{Fangfang Zhang}
\author[2]{Yaojun Chen}
\author[3]{Ervin Győri}
\author[4]{Xiutao Zhu}
\date{}

\affil[1]{School of Applied Mathematics, Nanjing University of Finance and Economics, Nanjing  210023, P.R. CHINA}
\affil[2]{Department of Mathematics, Nanjing University.}
\affil[3]{Alfr\'ed R\'enyi Institute of Mathematics, Hungarian Academy of Sciences. }
\affil[4]{School of mathematics, Nanjing University of Aeronautics and Astronautics, Nanjing  211106, P.R. CHINA}

 
\date{}
\maketitle
\begin{abstract}
The generalized Tur\'an number $\ex(n,K_s,F)$ denotes the maximum number of copies of $K_s$ in an $n$-vertex $F$-free graph. Let $kF$ denote $k$ disjoint copies of $F$. Gerbner, Methuku and Vizer [DM, 2019, 3130-3141] gave a lower bound for $\ex(n,K_3,2C_5)$ and obtained the magnitude of $\ex(n, K_s, kK_r)$. In this paper, we determine the exact value of $\ex(n,K_3,2C_5)$  and described the unique extremal graph for large $n$. Moreover, we also determine the exact value of $\ex(n,K_r,(k+1)K_r)$ which  generalizes some known results.
\vskip 2mm


\end{abstract}

\section{Introduction}
Let $G$ be a graph with the set of vertices $V(G)$. For two graphs $G$ and  $H$, let $G\cup H$ denote the disjoint union of $G$ and $H$, and $kG$ denote $k$ disjoint copies of $G$. We write $G+H$ for the join of $G$ and $H$, the graph obtained from $G\cup H$ by adding all edges between  $V(G)$ and $V(H)$. We use $K_n$, $C_n$, $P_n$ to denote the complete graph, cycle, and path on $n$ vertices, respectively. Let $K_s(G)$ denote the number of copies of $K_s$ in $G$.

For a graph $F$, the Tur\'an number of $F$, denote by $\ex(n,F)$, is the maximum number of edges in an $F$-free graph $G$ on  $n$ vertex. In 1941, Tur\'an \cite{turan} proved that the balanced complete $r$-partite graph on $n$ vertices, called Tur\'an graph $T_r(n)$, is the unique extremal graph of $\ex(n,K_{r+1})$. Starting from this, the Tur\'an problem has attracted a lot of attention. The study of  disjoint copies of a given graph in the context of Tur\'an  numbers is very rich. The first result is due to Erd\H{o}s and Gallai \cite{Gallai} who determined the Tur\'an number of $\ex(n,kK_2)$ for all $n$. Later Simonovits \cite{Simonovits} and independently Moon \cite{Moon} determined the Tur\'an number of disjoint copies of cliques. In \cite{Gorgol} Gorgol initiated the systematic investigation of Tur\'an numbers of disjoint copies of graphs and proved the following.
\begin{Theorem}(Gorgol \cite{Gorgol})\label{Thm1}
For every graph $F$ and $k\ge 1$, $$\ex(n,kF)=\ex(n,F)+O(n).$$
\end{Theorem}

In this paper we study the generalized Tur\'an number of disjoint copies of graphs. The generalized Tur\'an number $\ex(n,T,F)$ is the maximum number of copies of $T$ in any $F$-free graph on $n$ vertices. Obviously, $\ex(n,K_2,F)=\ex(n,F)$. The earliest result in this topic is due to Zykov \cite{Zykov} who proved that
$\ex(n,K_s,K_r)=K_s(T_{r-1}(n))$.
\begin{Theorem}(Zykov \cite{Zykov})\label{Zykov}
For all $n$,
$$\ex(n,K_s,K_r)=K_s(T_{r-1}(n)),$$
and $T_{r-1}(n)$ is the unique extremal graph.
\end{Theorem}
In recent years, the problem of estimating generalized Tur\'an number has received a lot of attention. Many classical results have been extended to generalized Tur\'an problem, see \cite{Alon, chase,Grzesik,Hatami,Luo,Ma,Wang, zhu}.

 Theorem \ref{Thm1} implies that the classical Tur\'an number $\ex(n,kF)$ and $\ex(n,F)$ always have the same order of magnitude. However, this is not true for generalized Tur\'an number.  The function $\ex(n,K_3,C_5)$ has attracted a lot of attentions, see \cite{Bollobas,Beka, Beka2}, the best known upper bound is given by Lv and Lu,
\begin{Theorem}(Lv and Lu \cite{Lv} )\label{Thm2}
$\ex(n,K_3,C_5)\le \frac{1}{2\sqrt{6}}n^{\frac{3}{2}}+o(n^{\frac{3}{2}})$.
\end{Theorem}
And  Gerbner, Methuku and Vizer \cite{Gerbner} proved $\ex(n,K_3,2C_5)=\Theta(n^{2})$ \cite{Gerbner}. This implies that the order of magnitudes of $\ex(n,H,F )$ and $\ex(n,H, kF)$ may differ.  They also obtained a lower bound for $\ex(n,K_3,2C_5)$ which is obtained by joining a vertex to a copy of $T_2(n-1)$.  In this paper, we show the graph $K_1+T_2(n-1)$  is indeed the unique extremal graph for $\ex(n,K_3,2C_5)$.

\begin{Theorem}\label{Thm3}
For sufficiently large $n$,
$$\ex(n,K_3,2C_5)=\left\lfloor\frac{(n-1)^2}{4}\right\rfloor,$$
and $K_1+T_2(n-1)$ is the unique extremal graph.
\end{Theorem}
%

We also focus on the generalized the Tur\'an number of disjoint copies of cliques. Since $\ex(n,K_s,K_r)$ is known \cite{Zykov}, it is natural to study the function $\ex(n,K_s,kK_r)$. Gerbner, Methuku and Vizer \cite{Gerbner} obtained the  asymptotic value of $\ex(n,K_s, kK_r)$.
\begin{Theorem}(Gerbner, Methuku and Vizer \cite{Gerbner})
If $s< r$, then
$$\ex(n,K_s,kK_r)=(1+o(1))\binom{r-1}{s}\left(\frac{n}{r-1}\right)^s.$$
If $s\ge r\ge 2$ and $k\ge 2$, then
$$\ex(n,K_s,kK_r)=\Theta(n^x),$$
where $x=\left\lceil\frac{kr-s}{k-1}\right\rceil-1$.
\end{Theorem}

Liu and Wang \cite{liu} determined the exact value of $\ex(n,K_r,2K_r)$ for $r\ge 3$ and $n$ sufficiently large. A new proof of  $\ex(n,K_r,2K_r)$ can be found in \cite{Yuan} by Yuan and Yang. Gerbner and Patk\'{o}s \cite{GP} determined $\ex(n,K_s,2K_r)$ for all $s\ge r\ge 3$ and $n$ sufficiently large. In this paper, we determine the value of $\ex(n,K_r,(k+1)K_r)$ for all $r\ge 2$, $k\ge 1$ and $n$ sufficiently large.
\begin{Theorem}\label{Thm5}
There exists a constant $n_0(k,r)$ depending on $k$ and $r\ge 2$ such that when $n\ge n_0(k,r)$,
$$\ex(n,K_r,(k+1)K_r)=K_r(K_k+T_{r-1}(n-k)),$$
and $K_k+T_{r-1}(n-k)$ is the unique extremal graph.
\end{Theorem}

 The detailed proofs of Theorems \ref{Thm3} and  \ref{Thm5} will be presented in Sections 3 and 4, respectively.


%

\section{Proof of Theorem \ref{Thm3}}
Suppose $n$ is large enough  and let $G$ be an $n$-vertex $2C_5$-free graph  with $\ex(n,K_3,2C_5)$ copies of  triangles. Since $K_1+T_2(n-1)$ contains no $2C_5$, thus $K_3(G)\ge \lfloor(n-1)^2/4\rfloor$. Next we will show that $G=K_1+T_2(n-1)$.  Since $n$ is sufficiently large and by Theorem \ref{Thm2}, $G$ must contain a copy of $C_5$, say $C=v_1v_2v_3v_4v_5v_1$. Then $G\setminus C$ contains no $C_5$. By Theorem \ref{Thm2} again, we have
$$K_3(G\setminus C)\le\frac{1}{2\sqrt{2}}(n-5)^{\frac{3}{2}}+o((n-5)^{\frac{3}{2}}).$$
We claim that there is at least one vertex in $V(C)$ whose neighborhood contains a copy of $6P_4$. To prove this, we need a theorem obtained by Bushaw and Kettle \cite{Kettle}.
\begin{Theorem}(Bushaw and Kettle\cite{Kettle})
For $k\ge 2$, $\ell\ge 4$ and $n\ge 2\ell+2k\ell(\lceil\ell/2\rceil+1)\binom{\ell}{\lfloor\ell/2\rfloor}$,
\[\ex(n,kP_\ell)=\binom{k\lfloor\ell/2\rfloor-1}{2}+(k\lfloor\ell/2\rfloor-1)(n-k\lfloor\ell/2\rfloor+1)+\lambda,\]
where $\lambda=1$ if $\ell$ is odd, and $\lambda=0$ if $\ell$ is even.
\end{Theorem}
By Theorem 7, we know $\ex(n,6P_4)\le \max\left\{\binom{872}{2},11(n-6)\right\}$. Now suppose no vertex in $V(C)$ contains $6P_4$ in its neighborhood. Then the number of triangles containing $v_i$ is at most
$$e(G[N(v_i)])\le \ex(n,6P_4)=11n+o(n).$$
Therefore, the total number of triangles satisfies
\[\begin{split}
K_3(G)&\le \frac{1}{2\sqrt{2}}n^{\frac{3}{2}}+o(n^{\frac{3}{2}})+55n+o(n)\\
&=\frac{1}{2\sqrt{2}}n^{\frac{3}{2}}+o(n^{\frac{3}{2}})\\
&<\frac{(n-1)^2}{4}.
\end{split}\]
The last inequality holds when $n$ is large. A contradiction.

Therefore, we may assume that $v_1$ is the vertex in $V(C)$ such that $G[N(v_1)]$ contains a copy of $6P_4$. If $G\setminus v_1$ contains a copy of $C_5$, then at least one copy of $P_4$ in $G[N(v_1)]$ does not intersect with this $C_5$ and hence we find two disjoint $C_5$, a contradiction. Thus $G\setminus v_1$ is $C_5$-free. So we have
\begin{equation}\label{2.1}
K_3(G)\le e(G\setminus v_1)+K_3(G\setminus v_1).
\end{equation}
So if we have $e(G\setminus v_1)+K_3(G\setminus v_1)\le  \left\lfloor\frac{(n-1)^2}{4}\right\rfloor$, then the proof is completed. To prove this,  we need the following lemma.
\begin{Lemma}\label{lemma}
Let $n\ge 2\binom{68}{3}$. If $G$ is a $C_5$-free graph on $n$ vertices, then
$$e(G)+K_3(G)\le \left\lfloor\frac{n^2}{4}\right\rfloor,$$
and equality holds if and only if $G=T_2(n)$.
\end{Lemma}
\pf For each integer $n$, let $G_n$ be a $C_5$-free graph of $n$ vertices such that $e(G_n)+K_3(G_n)$ is maximum. For every $n$, if $G_n$ is also triangle-free, then by Tur\'an Theorem \cite{turan}, $e(G_n)\le \left\lfloor\frac{n^2}{4}\right\rfloor$. Hence, $e(G_n)+K_3(G_n)\le \left\lfloor\frac{n^2}{4}\right\rfloor$ and equality holds if and only if $G_n=T_2(n)$, we are done.

Next we shall prove that from $n\ge 2\binom{68}{2}$, each $G_n$ is triangle-free. To do this, let us define a function
$$\phi(n):=e(G_n)+K_3(G_n)-\left\lfloor\frac{n^2}{4}\right\rfloor.$$
Since $T_2(n)$ is $C_5$-free and $e(T_2(n))+K_3(T_2(n))= \left\lfloor\frac{n^2}{4}\right\rfloor$, we have $\phi(n)\ge 0$. We claim that from $n\ge 68$, if $G_n$ contains a triangle, then 
\begin{align}\label{eq2.2}
\phi(n)<\phi(n-1)-1.    
\end{align}


 First suppose that $\delta(G_n)\ge \frac{n}{4}-1$. Let $xy$ be the edge of $G_n$ which is contained in the most number of triangles. Set $W=N(x)\cap N(y)=\{z_1,\ldots,z_w\}$. Since $G_n$ is $C_5$-free, $G_n[W]$ contains no edge unless $w\le 2$. Let $D_0=N(x)\setminus(W\cup\{y\})$, $D_i=N(z_i)\setminus(W\cup\{x,y\})$ for $1\le i\le w$  and $D_{w+1}=N(y)\setminus (W\cup \{x\})$.
We next show that $D_i$ satisfy the following properties for $0\le i\le w+1$.\medskip

\myitem{({\bf P1})} $|D_i|\ge \frac{n}{4}-w-2$ for $i=0,w+1$ and $|D_j|\ge \frac{n}{4}-4$ for $1\le j\le w$;

 \myitem{({\bf P2})} $D_i\cap D_j=\emptyset$ for $0\le i\not=j\le w+1$;

  \myitem{({\bf P3})} There are no edges between $D_i, D_j$.
\medskip

Since $\delta(G_n)\ge \frac{n}{4}-1$,  {\bf(P1)} is clearly true. Since $G_n$ is $C_5$-free, it is easy to see that $D_i\cap D_j=\emptyset$ for $1\le i\not=j\le w$. Suppose $D_0\cap D_i\neq \emptyset$ or $D_{w+1}\cap D_i\neq \emptyset$ for some $1\le i\le w$, by symmetry, let  $v\in D_0\cap D_i$. Then by the choice of $xy$, we have $w\ge 2$.  For $1\le j\le w$ and $j\neq i$, $vz_iyz_jxv$ is a  copy of $C_5$, a contradiction.  Thus {\bf(P2)} holds. Suppose $uv$ is an edge with $u\in D_i,v\in D_j$, then $uz_iyz_jvu$ is a copy of $C_5$ if $i,j\in[1,w]$,  $uz_iyxvu$ or $uz_ixyvu$ is a copy of $C_5$ if $i\in[1,w]$ and $j\in \{0,w+1\}$, $uxz_1yvu$ is a copy of $C_5$ if $i=0, j=w+1$, a contradiction. This implies  {\bf(P3)} holds.
\medskip

  Let $N=V(G_n)-W\cup \{x,y\}-\cup_{i=0}^{w+1} D_i$. By  {\bf(P1)} and {\bf(P2)}, we have
$$n=|N|+\sum_{i=0}^{w+1}|D_i|+w+2\ge |N|+2(\frac{n}{4}-w-2)+w(\frac{n}{4}-4)+w+2,$$
which implies $w\le 2$, $|N|\le \frac{n}{4}+7 $  and $D_i\not= \emptyset$ when $n\ge 61$. 
By the choice of $xy$, each vertex of $D_i$ has at most two neighbors in $G_n[D_i]$ for $0\le i\le w+1$ since there is no edge in $3$ triangles. By {\bf(P3)} and $\delta(G_n)\ge \frac{n}{4}-1$, each vertex in $D_i$ has at least $\frac{n}{4}-4$ neighbors in $N$. Let  $v_0\in D_0$ and $v_1\in D_{w+1}$. Because $n\ge 68$, we can deduce that $2(\frac{n}{4}-4)>\frac{n}{4}+7 \ge  |N|$ and hence $N(v_0)\cap N(v_1)\cap N\neq \emptyset$.  Then $uv_0xyv_1u$ is a copy of $C_5$, where $u\in N(v_0)\cap N(v_1)\cap N$, a contradiction.
We are done if the minimum degree is at least $\frac{n}{4}-1$.

Therefore, there is one vertex $v$ in $G_n$ such that $d(v)< \frac{n}{4}-1$ when $n\ge 68$. Because $G_n$ is $C_5$-free, $G_n[N(v)]$ is the disjoint union of stars and triangles which implies $e(G_n[N(v)])\le d(v)$. If we delete $v$ from $G_n$, it will destroy  at most $d(v)$ triangles and delete $d(v)$ edges. Hence,
\begin{align*}
&\phi(n-1)-\phi(n)\\
=&\left\lfloor\frac{n^2}{4}\right\rfloor-\left\lfloor\frac{(n-1)^2}{4}\right\rfloor-\{(e(G_n)+K_3(G_n))-(e(G_{n-1})+K_3(G_{n-1}))\} \notag\\
\ge& \frac{2n-2}{4}-\{(e(G_n)+K_3(G_n))-(e(G_{n}-v)+K_3(G_{n}-v))\}\notag\\
\ge& \frac{2n-2}{4}-2d(v)> \frac{2n-2}{4}-2(\frac{n}{4}-1)>1.
\end{align*} 
Hence our claim(inequality \ref{eq2.2}) holds for $n\ge 68$.

Note that  for $n_0\ge 68$, if $G_{n_0}$ contains no triangle, then $\phi(n_0)=0$. Moreover, for every $n\ge n_0$, we have that $G_n$ contains no triangles, either. Otherwise, we can find an integer $n$ such that $G_n$ contains a triangle but $G_{n-1}$ is triangle-free.  But then $\phi(n)\le \phi(n-1)-1<0$ by inequality \ref{eq2.2}, which is contrary to $\phi(n)\ge 0$. 
 Now let $n_0$ be the first integer after $68$ such that $G_{n_0}$ is triangle-free. Then
\[0\le \phi(n_0)\le \phi(n_0-1)-1< \phi(68)-(n_0-68)\le \binom{68}{2}+\binom{68}{3}+68-n_0.\]
This implies $n_0\le 2\binom{68}{3}$. Thus $G_n$ must be triangle-free for $n\ge 2\binom{68}{3}\ge n_0$. So $e(G_n)+K_3(G_n)=e(G_n)=\lfloor n^2/4\rfloor$ and $G_n=T_2(n)$ by Tur\'an Theorem \cite{turan}. The proof of Lemma \ref{lemma} is completed.$\hfill\square$

\vskip 3mm
Combining equation (\ref{2.1}) and Lemma \ref{lemma}, we can see that when $n$ is large, $K_3(G)\le \left\lfloor\frac{(n-1)^2}{4}\right\rfloor$ and equality holds if and only if $G=K_1+T_2(n-1)$.  The proof of Theorem \ref{Thm3} is completed. $\hfill\blacksquare$

\section{Proof of Theorem \ref{Thm5}}
We prove it by induction on $r$ and in each case, we always assume $n\ge n_{0}(k,r)=$. The base case  $r=2$ is the celebrated Erd\H{o}s-Gallai Theorem \cite{Gallai}, which says that
$$\ex(n,K_2,(k+1)K_2)=\max\left\{\binom{2k+1}{2},(n-k)k+\binom{k}{2}\right\}.$$
As $n\ge n_0(k,2)$, we know $\ex(n,K_2,(k+1)K_2)=K_2(K_k+T_1(n-k))$.

Let $r\ge 3$ and suppose that the result holds for all $r'<r$. Next we consider the case $\ex(n,K_r,(k+1)K_r)$. Let $G$ be a $(k+1)K_r$-free graph on $n$ vertices with $\ex(n,K_r,(k+1)K_r)$ copies of $K_r$. We may assume that $G$ contains $k$ disjoint copies of $K_r$. Otherwise we can add some edges  into $G$ unit the resulting graph contains  $k$ disjoint  $K_r$. But at least one $K_r$ in these $k$ disjoint $K_r$ is new which implies that the number of $K_r$ is increased, a contradiction. 
Let
$$I=\{X_1,\ldots,X_k\}$$ be a set of $k$ disjoint $r$-cliques in $G$, where $X_i$ is a copy of $K_r$. Let $V(I)=\cup_{i=1}^kV(X_i)$ and $N=G\setminus V(I)$. Clearly, $N$ contains no $K_r$. We say a vertex $v$ in $I$ is joined to an $(r-1)$-clique in $N$ if $v$ is adjacent to all vertices of this $(r-1)$-clique.  For each $X_i$, $i\in [k]$, we have the following property.

\begin{Claim}
Each $X_i$ contains at most one vertex which is joined to at least $kr+1$ disjoint $(r-1)$-cliques in $N$.    
\end{Claim}
\pf If not, suppose $u_1,u_1'\in V(X_1)$ are both joined to $kr+1$ disjoint $(r-1)$-cliques. First we can find an $(r-1)$-clique joined to $u_1$ in $N$. Since $u'_1$ is also joined to at least $kr+1$ disjoint $(r-1)$-cliques in $N$, we can find another $(r-1)$-clique joined to $u'_1$ which does not intersect with the $(r-1)$-clique joined to $u$.   Together with  $\{X_2,\ldots, X_k\}$, we find a copy of $(k+1)K_r$, a contradiction. $\hfill\square$

By Claim 1,
let $A=\{X_1,\ldots,X_a\}$ be a subset of $I$ such that there exists a vertex in $X_i$, say $u_i$,  that is joined to at least $kr+1$ disjoint $(r-1)$-cliques in $N$ for each $i\in [a]$.  Let $U=\{u_1,\ldots,u_a\}$.

Since $N$ is $K_r$-free, each $K_r$ in $G$ must intersect with some vertices in $V(I)$.  Then all $r$-cliques can be divided into two classes: the set of cliques in which all vertices  are contained in $V(N)\cup U$, ant the set of cliques containing at least one vertex  in $V(I)\setminus U$. We simply use $K_r(U)$ and $K_r(\overline{U})$ to denote the number of copies of $K_r$ in these two classes, respectively.

Suppose a $K_r$ in the first class contains $s$  vertices in $U$ and $r-s$ vertices in $N$, the number of $K_r$'s of this type is at most $\binom{a}{s}K_{r-s}(N)$. Since $N$ is $K_r$-free and by Theorem \ref{Zykov}, which says $\ex(n,K_s,K_r)=K_s(T_{r-1}(n))$,  we have $K_{r-s}(N)\le K_{r-s}\left (T_{r-1}(n-kr)\right)\le \binom{r-1}{r-s}\left(\frac{n-kr}{r-1}\right)^{r-s}$. Then
\begin{align}\label{3.1}
K_r(U)&\le \sum_{s=1}^r\binom{a}{s}K_{r-s}(N)\nonumber\\
&\le a\left(\frac{n-kr}{r-1}\right)^{r-1}+\binom{a}{2}\binom{r-1}{r-2}\left(\frac{n-kr}{r-1}\right)^{r-2}+O(n^{r-3}).
\end{align}

Next we calculate the size of $K_r(\overline{U})$. Each vertex $v\in V(I)\setminus U$ is joined to at most $kr$ independent $(r-1)$-cliques in $N$. Hence the number of $K_r$  containing $v$ and  $r-1$ vertices of $N$ is at most
\[\begin{split}
K_{r-1}(G[N(v)\cap V(N)])&\le \ex(n-kr,K_{r-1}, (kr+1)\cdot K_{r-1})\\
&=K_{r-1}\left (K_{kr}+T_{r-2}(n-2kr)\right )\\
&\le(kr)\left(\frac{n-2kr}{r-2}\right)^{r-2},
\end{split}\]
the second equality comes from the induction hypothesis. Any other copies of $K_r$ in $K_r(\overline{U})$ contains at most $r-2$ vertices in $N$ and at least one vertex in $V(I)\setminus U$. So the number of such $r$-cliques is at most
$$
\sum_{s=2}^r\left(\binom{kr}{s}-\binom{a}{s}\right)K_{r-s}(N)
\le \left(\binom{kr}{2}-\binom{a}{2}\right)\binom{r-1}{r-2}\left(\frac{n-kr}{r-1}\right)^{r-2}+O(n^{r-3}).$$
Hence, 
\begin{equation}\label{3.2}
K_r(\overline{U})\le \left(kr+\left(\binom{kr}{2}-\binom{a}{2}\right)\binom{r-1}{r-2}\right)\left(\frac{n-kr}{r-1}\right)^{r-2}+O(n^{r-3}).
\end{equation}

\vskip 3mm
Therefore, by inequality (\ref{3.1}) and (\ref{3.2}), we have
\begin{align}\label{3.3}
K_r(G)\le a\left(\frac{n-kr}{r-1}\right)^{r-1}+\left(kr+\binom{kr}{2}\binom{r-1}{r-2}\right)\left(\frac{n-kr}{r-1}\right)^{r-2}+O(n^{r-3}).
\end{align}
On the other hand, since $K_k+T_{r-1}(n-k)$ is $(k+1)K_r$-free, we know that
\begin{align}\label{3.4}
K_r(G)\ge k\left(\frac{n-k}{r-1}\right)^{r-1}+O(n^{r-2}).
\end{align}
When $n$ is greater than some constant $n_0(k,r)$, inequalites  (\ref{3.3}) and (\ref{3.4}) hold mean $a=k$ and then $U=\{u_1,\ldots,u_k\}$.

Let $G'=G\setminus U$. We claim that $G'$ is also $K_r$-free. Suppose not, $G'$ contains a $r$-clique, denote by $X_0'$. Since each $u_i$ is joined to at least $kr+1$ independent copies of  $K_{r-1}$'s in $N$, at least $(k-1)r+1$ of whom are disjoint with $X_0'$ for each $i\in [k]$. Then we can find a $r$-clique $X_1'$ such that  $u_1\in X_1'$ and $V(X_1')\cap V(X_0')=\emptyset$. Next, we claim that we may find another $k$ independent $r$-cliques such that each is disjoint with $X_0'$. 
Suppose we have found pairwise disjoint $r$-cliques $X_1',\ldots,X_{i-1}'$ such that $u_j\in X_j'$ for $j\in [i-1] $ and $i\le k$. 
Then,   in $G'[N(u_i)]$, there are   at least $(k-1)r+1-(i-1)(r-1)\ge 1$ independent $(r-1)$-cliques which disjoint with $\{X_0',X_1',\ldots,X_{i-1}'\}$.
That is we can choose a $(r-1)$-clique and thus a $r$-clique $X_i'$ such $u_i\in X_i'$ and $X_0',X_1',\ldots,X_{i}'$ are pairwise disjoint.  The procedure can keep going until we find $k$ independent $r$-cliques $X_1',\ldots,X_{k}'$. Then  $X_0',X_1',\ldots,X_{k}'$ forms a $(k+1)K_r$, a contradiction.

Since $G'$ is $K_r$-free, by Zykov's Theorem,  $K_{r-i}(G')\le K_{r-i}(T_{r-1}(n-k))$ and the equality holds if and only if  $G'=T_{r-1}(n-k)$. Thus
\[K_r(K_k+T_{r-1}(n-k))\le K_r(G)\le \sum_{i=0}^{r}\binom{k}{i}K_{r-i}(G')=K_r(K_k+T_{r-1}(n-k)).\]
The condition of the equality holds means $G=K_k+T_{r-1}(n-k)$.  The proof of Theorem \ref{Thm5} is completed. $\hfill\blacksquare$
%
%
\section{Acknowledge}
The research of the authors Gy\H{o}ri is partially supported by the National Research, Development and Innovation Office NKFIH, grants  K132696, SNN-135643 and K126853, Chen is supported by SNFC under grant numbers 12161141003 and 11931006, Zhang is supported by NSFC under grant numbers 12101298.


\end{document}